\input amstex
\documentstyle{amsppt}
%
%
\nopagenumbers
\def\Lj{L^{\raise 1pt\hbox{$\ssize j$}}}
\def\checkboldL{\kern -1pt\check{\kern 2pt\bold L}\kern -1pt}
\def\Ker{\operatorname{Ker}}
\def\rank{\operatorname{rank}}
\def\negskp{\hskip -2pt}
\pagewidth{360pt}
\pageheight{606pt}
\leftheadtext{Ruslan A. Sharipov}
\rightheadtext{On the subset of normality equations \dots}
\topmatter
\title On the subset of normality equations\\
\lowercase{describing generalized \uppercase{L}egendre
transformation.}
\endtitle
\author
R.~A.~Sharipov
\endauthor
\abstract
Normality equations describe Newtonian dynamical systems admitting
normal shift of hypersurfaces. They were first derived in Euclidean
geometry, then in Riemannian geometry. Recently they were rederived
in more general case, when geometry of manifold is given by
generalized Legendre transformation. As appears, in this case some
part of normality equations describe generalized Legendre transformation
itself irrespective to that Newtonian dynamical system, for which
others are written. In present paper this smaller part of normality
equations is studied.
\endabstract
\address Rabochaya street 5, 450003, Ufa, Russia
\endaddress
\email \vtop to 20pt{\hsize=280pt\noindent
R\_\hskip 1pt Sharipov\@ic.bashedu.ru\newline
r-sharipov\@mail.ru\vss}
\endemail
\urladdr
http:/\negskp/www.geocities.com/r-sharipov
\endurladdr
\endtopmatter
\loadbold
\TagsOnRight
\document
\head
1. Newtonian dynamical systems\\
and generalized Legendre transformation.
\endhead
    Let $M$ be smooth manifold of dimension $n$. We say that
the motion of a point $p=p(t)$ of this manifold obeys Newton's
second low if in local chart it is described by the following
ordinary differential equations:
$$
\xalignat 2
&\hskip -2em
\dot x^i=v^i,
&&\dot v^i=\Phi^i(x^1,\ldots,x^n,v^1,\ldots,v^n).
\tag1.1
\endxalignat
$$
Here $v^1,\,\ldots,\,v^n$ are components of velocity vector $\bold v$
of moving point. Its mass is assumed to be equal to unity: $m=1$.
Therefore functions $\Phi^1,\,\ldots,\,\Phi^n$ in \thetag{1.1} play
the role of force vector, though, unlike $v^1,\,\ldots,\,v^n$, they are
not components of tangent vector to $M$.\par
    Not always, but very often differential equations \thetag{1.1}
are associated with some extremal principle and hence are given 
implicitly by Euler-Lagrange equations:
$$
\xalignat 2
&\dot x^i=v^i,
&&\frac{d}{dt}\!\left(\frac{\partial L}
{\partial v^i}\right)=\frac{\partial L}{\partial x^i}
\endxalignat
$$
In this case they can be transformed to Hamiltonian form
$$
\xalignat 2
&\dot x^i=\frac{\partial H}{\partial p_i},
&&\dot p_i=-\frac{\partial H}{\partial x^i}
\endxalignat
$$
by means of classical Legendre transformation that relates velocity
vector $\bold v$ and momentum covector $\bold p$ according to the
following formula:
$$
\hskip -2em
p_i=\frac{\partial L}{\partial v^i}.
\tag1.2
$$
In \cite{1} and \cite{2} more general transformation was considered.
It is given by functions
$$
\hskip -2em
\cases p_1=L_1(x^1,\ldots,x^n,v^1,\ldots,v^n),\\
.\ .\ .\ .\ .\ .\ .\ .\ .\ .\ .\ .\ .\ .\ .\ .\
.\ .\ .\ .\ .\ .\ .\ \\
p_n=L_n(x^1,\ldots,x^n,v^1,\ldots,v^n).
\endcases
\tag1.3
$$
From geometric point of view generalized Legendre transformation
\thetag{1.3} is a smooth {\bf fiber-preserving} map from tangent
bundle to cotangent bundle:
$$
\hskip -2em
\lambda\!: TM\to T^*\!M.
\tag1.4
$$
Fiber-preserving means that each fixed fiber of tangent
bundle $TM$ is mapped into a fiber of $T^*\!M$ over the same base
point of $M$. For the sake of simplicity we shall assume generalized
Legendre map \thetag{1.4} to be diffeomorphic. Then inverse map
$$
\hskip -2em
\lambda^{-1}\!: T^*\!M\to TM
\tag1.5
$$
is also fiber-preserving. In local chart it is given by functions
$$
\hskip -2em
\cases v^1=V^1(x^1,\ldots,x^n,p_1,\ldots,p_n),\\
.\ .\ .\ .\ .\ .\ .\ .\ .\ .\ .\ .\ .\ .\ .\ .\
.\ .\ .\ .\ .\ .\ .\ \\
v^n=V^n(x^1,\ldots,x^n,p_1,\ldots,p_n).
\endcases
\tag1.6
$$\par 
     In paper \cite{1} generalized Legendre maps \thetag{1.4}
and \thetag{1.5} were used in order to transform dynamical system
\thetag{1.1} to $\bold p$-representation. Here it looks like
$$
\xalignat 2
&\hskip -2em
\dot x^i=V^i
&&\dot p_i=\Theta_i,
\tag1.7
\endxalignat
$$
where functions $V^1,\,\ldots,\,V^n$ are given by \thetag{1.6}, while
$\Theta_1,\,\ldots,\,\Theta_n$ are similar functions playing the same
role as function $\Phi^1,\,\ldots,\,\Phi^n$ in \thetag{1.1}.
Then in paper \cite{1} shift of hypersurfaces along trajectories of
dynamical system \thetag{1.7} was studied and {\bf theory of Newtonian
dynamical systems admitting normal shift of hypersurfaces} was
generalized to present non-metric geometry given by maps \thetag{1.4}
and \thetag{1.5}. Previous stage of development of this theory is
reflected in paper \cite{3} and in theses \cite{4} and \cite{5}
(see also recent papers \cite{6--13}).\par
    Main result of theory constructed in paper \cite{1} is a set
of {\bf normality equations}. This is rather huge system of partial
differential equations with respect to functions $V^1,\,\ldots,\,V^n$
and $\Theta_1,\,\ldots,\,\Theta_n$. In paper \cite{2} normality
equations were transformed back to $\bold v$-representation. Here
they form a system of partial differential equations with respect
to functions $\Phi^1,\,\ldots,\,\Phi^n$ in \thetag{1.1} and
functions $L_1,\,\ldots,\,L_n$ in \thetag{1.3}. Total set of normality
equations is divided into two parts: {\bf weak normality equations}
written for $n\geqslant 2$ and {\bf additional normality equations},
which are present only in multidimensional case $n\geqslant 3$.
Additional normality equations in turn are subdivided into {\bf three}
parts. It is remarkable that equations in the {\bf first part} have
no entries of functions $\Phi^1,\,\ldots,\,\Phi^n$ in them. They form
a system of partial differential equations with respect to functions
$L_1,\,\ldots,\,L_n$ that define generalized Legendre transformation
\thetag{1.4}. Further we shall call them {\bf normality equations for
generalized Legendre transformation}. Main goal of present paper is
to study these equations and describe generalized Legendre transformations
determined by their solutions.
\head
2. Normality equations\\
for generalized Legendre transformation.
\endhead
    Values of functions $L_1,\,\ldots,\,L_n$ in \thetag{1.3} form
components of covector $\bold p\in T^*_p(M)$ when their arguments
are fixed. However, they do not form components of traditional
covector field. They form so called {\bf extended covector field}.
\definition{Definition 2.1} Extended tensor field $\bold X$ of type
$(r,s)$ in $\bold v$-representation is a tensor-valued function
$\bold X=\bold X(q)$ with argument $q=(p,\bold v)$ in tangent bundle
$TM$ and with values in the following tensor space:
$$
T^r_s(p,M)=\overbrace{T_p(M)\otimes\ldots
\otimes T_p(M)}^{\text{$r$ times}}\otimes
\underbrace{T^*_p(M)\otimes\ldots
\otimes T^*_p(M)}_{\text{$s$ times}}.
$$
\enddefinition
\noindent
Extended covector field is a special case of extended tensor field,
when $r=0$ and $s=1$. Now we shall not discuss theory of extended
tensor fields, referring reader to Chapters~\uppercase
\expandafter{\romannumeral 2}, \uppercase\expandafter{\romannumeral
3}, and \uppercase\expandafter{\romannumeral 4} of thesis \cite{4}.
However, we should note that if
$$
X^{i_1\ldots\,i_r}_{j_1\ldots\,j_s}=X^{i_1\ldots\,i_r}_{j_1\ldots
\,j_s}(x^1,\ldots,x^n,v^1,\ldots,v^n)
$$
are components of extended tensor field $\bold X$, then partial
derivatives
$$
\hskip -2em
\tilde\nabla_{\!k}X^{i_1\ldots\,i_r}_{j_1\ldots\,j_s}=
\frac{\partial X^{i_1\ldots\,i_r}_{j_1\ldots\,j_s}}
{\partial v^k}
\tag2.1
$$
are components of another extended tensor field $\tilde\nabla\bold X$.
Therefore in \thetag{2.1} we use symbol of covariant derivative 
$\tilde\nabla_k$
for partial derivative $\partial/\partial v^k$.\par
    Let's apply covariant differentiation $\tilde\nabla$ to extended
covector field $\bold L$ with components \thetag{1.3}. As a result we
get extended tensor field $\bold g$ of type $(0,2)$ with components
$$
\hskip -2em
g_{qk}=\tilde\nabla_{\!k}L_q.
\tag2.2
$$
Matrix $g_{qk}$ in \thetag{2.2} is non-degenerate since it coincides with
Jacobi matrix for diffeomorphic map \thetag{1.4}. Hence we can consider 
inverse matrix with components $g^{qk}$. It defines extended tensor field 
of type $(2,0)$, we denote it by the same symbol $\bold g$. Though being
non-symmetric, tensor field $\bold g$ with components \thetag{2.2} and its
dual field with components $g^{qk}$ here play the same role as metric tensor 
and dual metric tensor in Riemannian geometry.\par
    Now, according to paper \cite{2}, we define extended scalar field $\Omega$
and operator-valued extended tensor field $\bold P$. They are determined as
follows:
$$
\xalignat 2
&\hskip -2em
\Omega=\sum^n_{s=1}L_s\,L^s=|\bold L|^2,
&&P^i_j=\delta^i_j-\frac{L^i\,L_j}{|\bold L|^2}.
\tag2.3
\endxalignat
$$
Here $L^s$ and $L^i$ are components of extended vector field $\bold L$
dual to covector field $\bold L$ with components \thetag{1.3} with
respect to non-symmetric metric \thetag{2.2}:
$$
\pagebreak
\hskip -2em
L^i=\sum^n_{s=1}L_s\,g^{si}.
\tag2.4
$$
Being more accurate, we should say that \thetag{2.4} are components of vector 
field
right-dual to covector field $\bold L$. One can also define left-dual vector 
field
with components 
$$
\hskip -2em
\check L^i=\sum^n_{s=1}g^{is}\,L_s.
\tag2.5
$$
In \thetag{2.3} we denoted $\Omega=|\bold L|^2$, it is positive
if non-symmetric metric \thetag{2.2} is positive. However, this
is not obligatory. We shall only require that $\Omega\neq 0$ since
it is in denominator in second formula \thetag{2.3}.\par
     Now we are ready to write normality equations for generalized Legendre 
transformation \thetag{1.4}. In local chart they are written as follows:
$$
\hskip -2em
\sum^n_{r=1}\sum^n_{s=1}(A^{rs}-A^{sr})\,P^i_r\,P^j_s=0.
\tag2.6
$$
Here $A^{rs}$ are components of extended tensor field $\bold A$. According to
paper \cite{2}, in $\bold v$-representation they are given by formula
$$
\hskip -2em
A^{rs}=\sum^n_{q=1}g^{qr}\,\tilde\nabla_{\!q}L^s.
\tag2.7
$$
Note that metric tensor $\bold g$ in \thetag{2.2}, projector field $\bold P$,
and tensor field $\bold A$ in \thetag{2.7} are completely determined by
covector field $\bold L$. Therefore \thetag{2.6} form a system of partial
differential equations with respect to functions $L_1,\,\ldots,\,L_n$.
Further steps are intended to study these equations. Note also that
equations \thetag{2.6} are written only for multidimensional case
$n\geqslant 3$. In two-dimensional case $n=2$ we have no restrictions
for generalized Legendre transformation \thetag{1.4}.
\head
3. Preliminary transformation of normality equations.
\endhead
    Let's consider formula \thetag{2.7}. Applying formula \thetag{2.4}
to $L^s$ in it, we derive the following expression for covariant derivative
$\tilde\nabla_{\!q}L^s$:
$$
\gathered
\tilde\nabla_{\!q}L^s=\tilde\nabla_{\!q}\left(\,\shave{\sum^n_{i=1}}
L_i\,g^{is}\!\right)=\sum^n_{i=1}\tilde\nabla_{\!q}L_i\,g^{is}+
\sum^n_{i=1}L_i\,\tilde\nabla_{\!q}g^{is}=\sum^n_{i=1}g_{iq}\,g^{is}
\,-\\
-\sum^n_{i=1}\sum^n_{a=1}\sum^n_{k=1}L_i\,g^{ia}\,\tilde\nabla_{\!q}
g_{ak}\,g^{ks}=\sum^n_{i=1}g_{iq}\,g^{is}-\sum^n_{i=1}\sum^n_{a=1}
\sum^n_{k=1}L_i\,g^{ia}\,\tilde\nabla_{\!q}\!\tilde\nabla_{\!k}L_a
\,g^{ks}.
\endgathered
$$
Upon substituting this expression into \thetag{2.7} for $A^{rs}$ we
obtain
$$
A^{rs}=g^{rs}-\sum^n_{a=1}\sum^n_{q=1}\sum^n_{k=1}g^{qr}\,g^{ks}\,
L^a\,\tilde\nabla_{\!q}\!\tilde\nabla_{\!k}L_a.
\tag3.1
$$
It is obvious that last term in \thetag{3.1} is symmetric with respect
to indices $r$ and $s$. Therefore it makes no contribution to ultimate
form of normality equations when we substitute \thetag{3.1} into
\thetag{2.6}. Thus from \thetag{2.6} we derive
$$
\hskip -2em
\sum^n_{r=1}\sum^n_{s=1}(g^{rs}-g^{sr})\,P^i_r\,P^j_s=0.
\tag3.2
$$
Now let's apply formula \thetag{2.3} to components of projector field
$P^i_r$ and $P^j_s$ in \thetag{3.2}:
$$
\gather
\sum^n_{r=1}\sum^n_{s=1}(g^{rs}-g^{sr})\,P^i_r\,P^j_s=
\sum^n_{r=1}\sum^n_{s=1}(g^{rs}-g^{sr})\left(\delta^i_r-
\frac{L^i\,L_r}{|\bold L|^2}\right)\!\left(\delta^j_s-
\frac{\Lj\,L_s}{|\bold L|^2}\right)=\\
\vspace{1ex}
=g^{ij}-g^{ji}-\frac{(\check L^i-L^i)\,L^j}{|\bold L|^2}
-\frac{(\Lj-\check\Lj)\,L^i}{|\bold L|^2}=g^{ij}-g^{ji}
-\frac{\check L^i\,\Lj-L^i\,\check\Lj}{|\bold L|^2}=0.
\endgather
$$
Here $L^i$, $\Lj$, $\check L^i$, and $\check\Lj$ are determined by
formulas \thetag{2.4} and \thetag{2.5}. As a result of the above
calculations normality equations \thetag{2.6} are written as
$$
\hskip -2em
g^{ij}-\frac{\check L^i\,\Lj}{|\bold L|^2}=g^{ji}
-\frac{\check\Lj\,L^i}{|\bold L|^2}.
\tag3.3
$$
If we denote by $u^{ij}$ left hand side of the equality \thetag{3.3},
then $g^{ij}$ is given by formula
$$
\hskip -2em
g^{ij}=u^{ij}+\frac{\check L^i\,\Lj}{|\bold L|^2},
\tag3.4
$$
while normality equations \thetag{3.3} themselves are equivalent to
symmetry of tensor $\bold u$ with components $u^{ij}$. Thus,
non-symmetric metric $\bold g$ is expressed through symmetric
tensor $\bold u$ by formula \thetag{3.4}. This is basic observation
for the next step.
\head
4. Fine structure of metric tensor.
\endhead
    Let's fix some point $q=(p,\bold v)$ of $TM$ such that $|\bold L|\neq
0$. This means that we fix arguments of extended tensor
fields in \thetag{3.4}. Then values of $\bold g$ and $\bold u$ for that
fixed argument $q$ are tensors from $T^2_0(p,M)$, while values of
$\bold L$ and $\checkboldL$ are vectors from tangent space $T_p(M)$.
Tensors $\bold g$ and $\bold u$ of type $(2,0)$ can be treated as bilinear
forms (bilinear functions) with arguments in cotangent space $T^*_p(M)$:
$$
\xalignat 2
&\hskip -2em
\bold g=\bold g(\bold x,\bold y),
&&\bold u=\bold u(\bold x,\bold y).
\tag4.1
\endxalignat
$$
Due to symmetry $u^{ij}=u^{ji}$ bilinear form $\bold u$ in \thetag{4.1}
is symmetric, i\.\,e\.
$$
\bold u(\bold x,\bold y)=\bold u(\bold y,\bold x).
$$
It is well known fact from linear algebra (see \cite{14}) that each
symmetric bilinear form can be diagonalized. This means that one can
choose some special base in $T_p(M)$ and its dual base in $T^*_p(M)$
such that matrix $u^{ij}$ is diagonal
$$
\hskip -2em
u^{ij}=\Vmatrix \varepsilon_1 & 0 &\hdots & 0\\
0 & \varepsilon_2 &\hdots & 0\\
\vdots & \vdots &\ddots & \vdots\\
0 & 0 &\hdots & \varepsilon_n
\endVmatrix.
\tag4.2
$$
Here it is important to note that tensor field $\bold u$ is diagonalized
at one fixed point $q=(p,\bold v)$, not in whole neighborhood of that
point.
\proclaim{Lemma 4.1} For each point $q\in TM$, where $|\bold L|\neq 0$,
tensor field $\bold u$ and its matrix \thetag{4.2} are degenerate, i\.\,
e\. at least one number among $\varepsilon_1,\,\ldots,\,\varepsilon_n$
is equal to zero.
\endproclaim
\demo{Proof}Let's multiply both sides of \thetag{3.4} by $L_j$ and then
sum up with respect to index $j$. As a result we get the following
equality:
$$
\hskip -2em
\check L^i=\sum^n_{j=1}g^{ij}\,L_j=\sum^n_{j=1}\left(\!u^{ij}+
\frac{\check L^i\,\Lj}{|\bold L|^2}\right)\!L_j=\sum^n_{j=1}u^{ij}\,L_j+
\check L^i.
\tag4.3
$$
Comparing left and right hand sides of the equality \thetag{4.3},
we derive
$$
\hskip -2em
\sum^n_{j=1}u^{ij}\,L_j=0.
\tag4.4
$$
If $\bold L\neq 0$, then the equality \thetag{4.4} means that
$\det\bold u=0$. This proves lemma for all points $q=(p,\bold v)$,
where $\bold L\neq 0$. But $|\bold L|\neq 0$ implies that covector
$\bold L$ is non-zero. Thus, lemma~4.1 is proved.\qed\enddemo
\remark{Remark} Normality equation \thetag{2.6} is derived only for
those points, where $|\bold L|\neq 0$ (see \cite{1} and \cite{2}).
Indeed, $|\bold L|$ is in denominator in formula \thetag{2.3} for
$P^i_j$. Hence lemma~4.1 is sufficient result for our further
purposes.
\endremark
    Lemma~4.1 means that bilinear form $\bold u$ has nonzero kernel.
This is linear subspace in cotangent space $T^*_p(M)$ defined as
follows:
$$
\hskip -2em
\Ker\bold u=\{\bold x\in T^*_p(M)\!:\,\bold u(\bold x,\bold y)=0
\ \forall\bold y\in T^*_p(M)\}.
\tag4.5
$$
In terms of kernel \thetag{4.5} the equality \thetag{4.4} now can
be written as
$$
\hskip -2em
\bold L\in\Ker\bold u\neq \{0\}.
\tag4.6
$$
\proclaim{Lemma 4.2} For symmetric bilinear form $\bold u$ in
$T^*_p(M)$ defined by \thetag{3.4} its rank is $n-1$ and the
dimension of its kernel is equal to unity, i\.\,e\.
$$
\xalignat 2
&\hskip -2em
\rank\bold u=n-1,
&&\dim\Ker\bold u=1.
\tag4.7
\endxalignat
$$
\endproclaim
\demo{Proof} Let's multiply both sides of \thetag{3.4} by $P^s_j$
and sum up with respect to double index $j$. As a result we obtain
the following equality
$$
\hskip -2em
\sum^n_{j=1}P^s_j\,g^{ij}=\sum^n_{j=1}u^{ij}\,P^s_j=
\sum^n_{j=1}u^{ij}\left(\delta^s_j-\frac{L_j\,L^s}
{|\bold L|^2}\right)=u^{is}.
\tag4.8
$$
Here in the above calculations we used \thetag{4.4}. Sum in left hand
side of \thetag{4.8} represents matrix product of two matrices: $P^s_j$
and $g^{ij}$ transposed. Matrix $g^{ij}$ is non-degenerate, while rank
of projection operator $\bold P$ is equal to $n-1$. This proves the
equalities \thetag{4.7} and lemma~4.2 in whole.
\qed\enddemo
\nopagebreak
\proclaim{Lemma 4.3} Matrix equality $g^{ij}=u^{ij}+A^i\,\Lj$ is
equivalent to normality equations \thetag{3.2} if and only if matrix
$u^{ij}$ is symmetric and degenerate.
\endproclaim
\pagebreak
\demo{Proof} Above we have derived the equality \thetag{3.4} from
normality equation \thetag{3.2} and we have proved that matrix
$u^{ij}$ in \thetag{3.4} is degenerate (see lemma~4.1 and lemma~4.2).
Denoting $A^i=\check L^i/|\bold L|^2$ we get the equality $g^{ij}=u^{ij}+A^i
\,\Lj$. Thus, direct proposition of lemma~4.3 is proved.\par
     Let's prove converse proposition. Suppose that metric tensor
is given by the equality $g^{ij}=u^{ij}+A^i\,\Lj$, where $\Lj$ are
determined by formula \thetag{2.4}, $A^i$ are components of some vector,
while matrix $u^{ij}$ is symmetric and degenerate. Then there exists
some covector $\bold x\neq 0$ with components $x_1,\,\ldots,\,x_n$ such
that
$$
\xalignat 2
&\hskip -2em
\sum^n_{j=1}u^{ij}\,x_j=0,
&&\sum^n_{i=1}x_i\,u^{ij}=0.
\tag4.9
\endxalignat
$$
Applying relationships \thetag{4.9} to the equality $g^{ij}=u^{ij}
+A^i\,\Lj$, we get
$$
\aligned
\hskip -2em
x^j=\sum^n_{i=1}x_i\,g^{ij}=\sum^n_{i=1}x_i\,A^i\,\Lj=
\left<\bold x\,|\,\bold A\right>\cdot\Lj,\\
\check x^i=\sum^n_{j=1}g^{ij}\,x_j=\sum^n_{j=1}A^i\,\Lj\,x_j=
\left<\bold x\,|\,\bold L\right>\cdot A^i.
\endaligned
\tag4.10
$$
From first equality \thetag{4.10} we derive that covectors $\bold x$
and $\bold L$ are collinear:
$$
\hskip -2em
x_i=\sum^n_{j=1}x^j\,g_{ji}=\sum^n_{j=1}\left<\bold x\,|\,\bold A
\right>\,\Lj\,g_{ji}=\left<\bold x\,|\,\bold A\right>\cdot L_i.
\tag4.11
$$
Note that $\bold x\neq 0$ and $\bold L\neq 0$. Hence $\left<\bold x\,|
\,\bold A\right>\neq 0$. Substituting formula \thetag{4.11} for $x_j$
into both sides of second equality \thetag{4.10}, we obtain
$$
\hskip -2em
\left<\bold x\,|\,\bold A\right>\cdot\check L^i=
\left<\bold x\,|\,\bold A\right>\cdot|\bold L|^2\cdot A^i.
\tag4.12
$$
Since $\left<\bold x\,|\,\bold A\right>\neq 0$, we can cancel this factor
in \thetag{4.12}. Then we get formula for $A^i$:
$$
\hskip -2em
A^i=\frac{\check L^i}{|\bold L|^2}.
\tag4.13
$$
Substituting \thetag{4.13} back into the equality $g^{ij}=u^{ij}+A^i
\,\Lj$, we get formula coinciding with \thetag{3.4}. Using symmetry
of $u^{ij}$, we can transform it to \thetag{3.3}. Then multiplying
\thetag{3.3} by $P^r_i\,P^s_j$, upon summation with respect to
double indices $r$ and $s$ we rederive normality equations \thetag{3.2}.
Lemma~4.3 is proved.
\qed\enddemo
    Now let's multiply \thetag{3.4} by $g_{ir}\,g_{js}$ and let's sum
resulting equality with respect to double indices $i$ and $j$. Then we
introduce the following notations:
$$
\xalignat 2
&\hskip -2em
u_{sr}=\sum^n_{i=1}\sum^n_{j=1}u^{ij}\,g_{ir}\,g_{js},
&&\check L_r=\sum^n_{i=1}\check L^i\,g_{ir}.
\tag4.14
\endxalignat
$$
Here $\check L_1,\,\ldots,\,\check L_n$ are components of extended
\pagebreak
covector field left dual to vector field $\checkboldL$\,, while vector
field $\checkboldL$ is right dual to initial covector field $\bold L$.
In terms of these newly introduced notations \thetag{4.14} transformed
equality \thetag{3.4} is written as
$$
\hskip -2em
g_{sr}=u_{sr}+\frac{L_s\,\check L_r}{|\bold L|^2}.
\tag4.15
$$
Matrix $u_{sr}$ in \thetag{4.15} is symmetric. This matrix is
degenerate, its rank is equal to $n-1$. This follows from
\thetag{4.6} due to \thetag{4.14}. Moreover, $g_{ir}$ and $g_{js}$
in \thetag{4.14} are components of non-degenerate matrix, therefore
\thetag{4.15} is equivalent to \thetag{3.4}.
\proclaim{Lemma 4.4} Matrix equality $g_{sr}=u_{sr}+L_s\,A_r$ is
equivalent to normality equations \thetag{3.2} if and only if matrix
$u_{sr}$ is symmetric and degenerate. 
\endproclaim
\demo{Proof} Note that matrix equality \thetag{4.15} with symmetric
degenerate matrix $u_{sr}$, which was derived above from normality
equation \thetag{3.4}, is particular form of the equality $g_{sr}
=u_{sr}+L_s\,A_r$, where $A_r=\check L_r/|\bold L|^2$. This means
that direct proposition of lemma~4.4 is proved.\par
     Let's prove converse proposition. Suppose that metric tensor
is given by the equality $g_{sr}=u_{sr}+L_s\,A_r$, where matrix
$u^{ij}$ is symmetric and degenerate. Then there exists some vector
$\bold X\neq 0$ with components $X^1,\,\ldots,\,X^n$ such that
$$
\xalignat 2
&\hskip -2em
\sum^n_{r=1}u_{sr}\,X^r=0,
&&\sum^n_{s=1}X^s\,u_{sr}=0.
\tag4.16
\endxalignat
$$
Applying relationships \thetag{4.16} to the equality $g_{sr}=u_{sr}
+L_s\,A_r$, we get
$$
\aligned
\hskip -2em
\check X_s=\sum^n_{r=1}g_{sr}\,X^r=\sum^n_{r=1}L_s\,A_r\,X^r=
L_s\cdot \left<\bold A\,|\,\bold X\right>,\\
X_r=\sum^n_{s=1}X^s\,g_{sr}=\sum^n_{s=1}X^s\,L_s\,A_r=
A_r\cdot\left<\bold L\,|\,\bold X\right>.
\endaligned
\tag4.17
$$
From first equality \thetag{4.17} we derive that vectors $\bold X$
and $\checkboldL$ are collinear:
$$
\hskip -2em
X^r=\sum^n_{s=1}g^{rs}\,\check X_s=\sum^n_{s=1}\left<\bold A\,|\,
\bold X\right>\,g^{rs}\,L_s=\left<\bold A\,|\,\bold X\right>
\cdot\check L^r.
\tag4.18
$$
Note that $\bold X\neq 0$ and $\checkboldL\neq 0$. Hence $\left<\bold A
\,|\,\bold X\right>\neq 0$. Substituting formula \thetag{4.18} for $X^s$
into both sides of second equality \thetag{4.17} and taking into
account \thetag{4.14}, we get
$$
\hskip -2em
\left<\bold A\,|\,\bold X\right>\cdot\check L_r=
\left<\bold A\,|\,\bold X\right>\cdot|\bold L|^2\cdot A_r.
\tag4.19
$$
Since $\left<\bold A\,|\,\bold X\right>\neq 0$, we can cancel this factor
in \thetag{4.19}. As a result we obtain
$$
\hskip -2em
A_r=\frac{\check L_r}{|\bold L|^2}.
\tag4.20
$$
Substituting \thetag{4.20} back into the equality $g_{sr}=u_{sr}+L_s
\,A_r$, \pagebreak we get formula coinciding with \thetag{4.15}.
Remember that \thetag{4.15} is equivalent to \thetag{3.4} (see above). 
Further from \thetag{3.4} we can rederive normality equations
\thetag{3.2}. This step is the same as in proving previous lemma~4.3.
Thus, lemma~4.4 is proved.\qed\enddemo
\head
5. Skew symmetry and differential forms.
\endhead
    Now we shall draw some conclusions from lemma~4.4. Lemma~4.4
asserts that functions $L_1,\,\ldots,\,L_n$ of the form \thetag{1.3}
define generalized Legendre transformation \thetag{1.4} satisfying
normality equations \thetag{2.6} if and only if their partial
derivatives $g_{sr}=\tilde \nabla_{\!r}L_s$ are related to them by
means of the equality
$$
\hskip -2em
\frac{\partial L_s}{\partial v^r}=u_{sr}+L_s\,A_r,
\tag5.1
$$
where $u_{sr}$ are components of some symmetric degenerate extended
tensor field $\bold u$, which is not initially predefined, and $A_r$
are components of some extended covector field $\bold A$, which also
is not initially predefined. Alternating \thetag{5.1}, we get
$$
\hskip -2em
\frac{\partial L_s}{\partial v^r}-\frac{\partial L_r}{\partial v^s}=
L_s\,A_r-L_r\,A_s.
\tag5.2
$$
For matrix $u_{sr}$ due to its symmetry $u_{sr}=u_{rs}$ from
\thetag{5.1} we derive
$$
\hskip -2em
u_{sr}=\frac{1}{2}\left(\frac{\partial L_s}{\partial v^r}
+\frac{\partial L_s}{\partial v^r}\right)-
\frac{L_s\,A_r+L_r\,A_s}{2}.
\tag5.3
$$
If functions $A_1,\,\ldots,\,A_n$ are given, then \thetag{5.2} can
be treated as differential equations for functions $L_1,\,\ldots,\,L_n$.
Suppose we take some covector field $\bold A$ and solve differential
equations \thetag{5.2}. Does it mean that we can reconstruct the equality
\thetag{5.1} and further get the solution of normality equations
\thetag{2.6}\,? Indeed, we could define matrix $u_{rs}$ by formula
\thetag{5.3} and then derive \thetag{5.1} from \thetag{5.2} and
\thetag{5.3}. Anyway, matrix $u_{rs}$ determined by formula \thetag{5.3}
is symmetric, but it could be non-degenerate. In this case lemma~4.4
is not applicable and further thread of reasoning is torn.\par
    However, thing are not so bad. Note that partial differential equations
\thetag{5.2} admit gauge transformations of the following form:
$$
\xalignat 2
&\hskip -2em
L_r\to L_r,
&&A_r\to A_r-\lambda\,L_r.
\tag5.4
\endxalignat
$$
Here $\lambda$ is some scalar factor, i\.\,e\. some extended scalar field
in $M$. Applying gauge transformation \thetag{5.4} we get new fields
$\bold A'$ and $\bold u'$ from initial ones:
$$
\xalignat 2
&\hskip -2em
A'_r=A_r-\lambda\,L_r.
&&u'_{sr}=u_{sr}+\lambda\,L_s\,L_r.
\tag5.5
\endxalignat
$$
If matrix $u_{sr}$ in \thetag{5.5} is non-degenerate, then we can
calculate determinant of $u'_{sr}$:
$$
\hskip -2em
\det(u'_{sr})=\det(u'_{sr})\left(\!1+\lambda\shave{\sum^n_{s=1}}
\shave{\sum^n_{r=1}}w^{rs}\,L_r\,L_s\!\right)=0.
\tag5.6
$$
Here $w^{rs}$ is inverse matrix for $u_{rs}$. \pagebreak
Looking at characteristic
equation \thetag{5.6}, we see that it is linear with respect to scalar
factor $\lambda$. This means that it is solvable if and only if double
sum in round brackets is nonzero:
$$
\hskip -2em
\Vert\bold L\Vert_{\bold u}=\sum^n_{s=1}\sum^n_{r=1}w^{rs}\,
L_r\,L_s\neq 0.
\tag5.7
$$
Now we shall leave inequality \thetag{5.7} for separate study in
separate paper and we shall formulate main result of this section
in the following theorem.
\proclaim{Theorem 5.1} Any solution of differential equations
\thetag{5.2} defines locally diffeomorphic generalized Legendre
map \thetag{1.4} if metric tensor \thetag{2.2} is non-degenerate
and if one of the following two conditions is fulfilled: matrix
\thetag{5.3} is degenerate or $\Vert\bold L\Vert_{\bold u}\neq 0$,
if matrix \thetag{5.3} is non-degenerate.
\endproclaim
     Note that differential equations \thetag{5.2} have no partial
derivatives with respect to $x^1,\,\ldots,\,x^n$. This means that
we can fix some arbitrary point $p\in M$ and consider partial
differential equations \thetag{5.2} within fixed fiber of tangent
bundle. Then extended covector fields $\bold L$ and $\bold A$ can
be treated as differential $1$-forms:
$$
\xalignat 2
&\hskip -2em
\bold L=\sum^n_{i=1}L_i\,dv^i,
&&\bold A=\sum^n_{i=1}A_i\,dv^i.
\tag5.8
\endxalignat
$$
In terms of differential forms \thetag{5.8} differential equations
\thetag{5.2} are written as
$$
\hskip -2em
d\,\bold L=\bold L\wedge\bold A.
\tag5.9
$$
\remark{Remark}Here we should especially emphasize that differential
forms \thetag{5.8} are defined only within separate fibers of tangent
bundle $TM$. They cannot be canonically extended as $1$-forms in $TM$
in whole.
\endremark
\head
6. Compatibility conditions.
\endhead
    Initial normality equations \thetag{2.6}, as well as their
transformed counterparts \thetag{5.9}, form overdetermined system
of partial differential equations for the functions \thetag{1.3}.
They should be studied for compatibility. Let's apply external
differentiation operator $d$ to both sides of \thetag{5.9}. As a
result we get
$$
0=d\,(d\,\bold L)=d\,\bold L\wedge\bold A-\bold L\wedge d\,\bold A=
\bold L\wedge\bold A\wedge\bold A-\bold L\wedge d\,\bold A=
-\bold L\wedge d\,\bold A.
$$
This means that external product $\bold L\wedge d\,\bold A$ is equal
to zero:
$$
\bold L\wedge d\,\bold A=0
\tag6.1
$$
\proclaim{Lemma 6.1} For $1$-form $\bold L\neq 0$ and differential
$m$-form $\Omega$ the equality $\bold L\wedge\Omega=0$ is equivalent
to the equality $\Omega=\bold L\wedge\bold B$ for some differential
$(m-1)$-form $\bold B$.
\endproclaim
     Lemma~6.1 is special case of division theorem by E.~Cartan, see
proof in \cite{15}. Applying lemma~6.1 to $\Omega=d\,\bold A$ in
\thetag{6.1}, we get the equality
$$
\hskip -2em
d\,\bold A=\bold L\wedge\bold B,
\tag6.2
$$
where $\bold B$ is some differential $1$-form within separate fibers of
tangent bundle $TM$. Differential equations \thetag{6.2} form compatibility
condition for equations \thetag{5.9}. They have almost the same shape as
\thetag{5.9}. Therefore we shall treat them similarly:
$$
0=d\,(d\,\bold A)=d\,\bold L\wedge\bold B-\bold L\wedge d\,\bold B=
\bold L\wedge\bold A\wedge\bold B-\bold L\wedge d\,\bold B=\bold L
\wedge(\bold A\wedge\bold B-d\,\bold B).
$$
Applying lemma~6.1 to the above equality, we get differential equations
for $\bold B$:
$$
\hskip -2em
d\,\bold B=\bold A\wedge\bold B+\bold L\wedge\bold C.
\tag6.3
$$
Here $\bold C$ is some other $1$-form. Differential equations
\thetag{6.3} are a little bit more complicated than \thetag{5.9}
and \thetag{6.1}. But nevertheless we apply operator $d$ to them:
$$
\gather
d\,(d\,\bold B)=d\,\bold A\wedge\bold B-\bold A\wedge d\,\bold B
+d\,\bold L\wedge\bold C-\bold L\wedge d\,\bold C
=\bold L\wedge\bold B\wedge\bold B-\\
-\bold A\wedge\bold A\wedge\bold B-\bold A\wedge\bold L\wedge
\bold C+\bold L\wedge\bold A\wedge\bold C-\bold L\wedge d\,\bold C.
\endgather
$$
Applying lemma~6.1 to the above equality, we get differential equations
for $\bold C$:
$$
\hskip -2em
d\,\bold C=2\,\bold A\wedge\bold C+\bold L\wedge\bold D.
\tag6.4
$$
Now again, we apply external differentiation $d$ to the equations
\thetag{6.4} and we get
$$
\gather
d\,(d\,\bold C)=2\,d\,\bold A\wedge\bold C
-2\,\bold A\wedge d\,\bold C
+d\,\bold L\wedge\bold D
-\bold L\wedge d\,\bold D=2\,\bold L\wedge\bold B\wedge\bold C-\\
-4\,\bold A\wedge\bold A\wedge\bold C
-2\,\bold A\wedge\bold L\wedge\bold D
+\bold L\wedge\bold A\wedge\bold D
-\bold L\wedge d\,\bold D.
\endgather
$$
Applying lemma~6.1 to this equality, we derive differential equations
for $\bold D$:
$$
\hskip -2em
d\,\bold D=3\,\bold A\wedge\bold D+2\,\bold B\wedge\bold C
+\bold L\wedge\bold E.
\tag6.5
$$
Now it is clear that further steps require special notations and
study of recurrent procedure underlying all above formulas \thetag{5.9},
\thetag{6.2}, \thetag{6.3}, \thetag{6.4}, \thetag{6.5}. Let's denote
$$
\xalignat 3
&\hskip -2em
\bold L=\bold A_0, &&\bold A=\bold A_1, &&\bold B=\bold A_2,\\
\vspace{-1.5ex}
\tag6.6\\
\vspace{-1.5ex}
&\hskip -2em
\bold C=\bold A_3, &&\bold D=\bold A_4, &&\bold E=\bold A_5.
\endxalignat
$$
In terms of notations \thetag{6.6} introduced just above we can rewrite
our equations as
$$
\xalignat 2
&\hskip -2em
d\,\bold A_0=\bold A_0\wedge \bold A_1,
&&d\,\bold A_1=\bold A_0\wedge \bold A_2,\\
\vspace{-1.5ex}
\tag6.7\\
\vspace{-1.5ex}
&\hskip -2em
d\,\bold A_2=\bold A_0\wedge \bold A_3+\bold A_1\wedge \bold A_2,
&&d\,\bold A_3=\bold A_0\wedge \bold A_4+2\,\bold A_1\wedge \bold A_3.
\endxalignat
$$
Equations \thetag{6.5} are a little bit more complicated. They are
written as follows:
$$
\hskip -2em
d\,\bold A_4=\bold A_0\wedge\bold A_5+3\,\bold A_1\wedge\bold A_4
+2\,\bold A_2\wedge\bold A_3.
\tag6.8
$$
Looking at \thetag{6.7} and \thetag{6.8}, one can formulate a conjecture
concerning general structure of all such \pagebreak equations, for those,
which are already written, and for all others.
\proclaim{Conjecture 6.1} Differential equations \thetag{5.9} lead to
infinite series of compatibility conditions that in terms of notations
\thetag{6.6} can be written as
$$
\hskip -2em
d\,\bold A_k=\sum^{\left[\!\frac{k}{2}\!\right]}_{i=0}C^i_{k+1}
\,\bold A_i\wedge\bold A_{k+1-i}\text{, \ where \ }C^1_{k+1}=1.
\tag6.9
$$
Here $C^i_{k+1}$ are some constants similar to binomial coefficients,
but not coinciding with them. They should be calculated recurrently.
By square brackets in upper limit of sum in \thetag{6.9} we denote
entire part of fraction $k/2$.
\endproclaim
     First of all let's derive recurrent relationships for coefficients
$C^i_{k+1}$ in \thetag{6.9}. Applying external differentiation $d$ to
both sides of \thetag{6.9}, we get
$$
\gathered
0=d\,(d\,\bold A_k)=\sum^{\left[\!\frac{k}{2}\!\right]}_{i=0}
C^i_{k+1}\left(d\,\bold A_i\wedge\bold A_{k+1-i}-\bold A_i\wedge
d\,\bold A_{k+1-i}\right)=\\
=\bold A_0\wedge\bold A_1\wedge\bold A_{k+1}-\bold A_0\wedge
d\,\bold A_{k+1}
+\sum^{\left[\!\frac{k}{2}\!\right]}_{i=1}C^i_{k+1}\left(C^0_{i+1}
\,\bold A_0\wedge\bold A_{i+1}\,+\right.\\
\left.\vphantom{C^0_{i+1}\bold A_0\wedge\bold A_{i+1}}
+\,\dots\right)\wedge\bold A_{k+1-i}
-\sum^{\left[\!\frac{k}{2}\!\right]}_{i=1}C^i_{k+1}\,\bold A_i\wedge
\left(C^0_{k+2-i}\,\bold A_0\wedge\bold A_{k+2-i}+\dots\right)\!.
\endgathered\quad
\tag6.10
$$
Terms denoted by dots in the above equality have no entry of $\bold A_0$.
Below we shall prove that they do cancel each other. Now from \thetag{6.10}
we derive
$$
\hskip -2em
\aligned
A_0\,\wedge\!&\left(-d\,\bold A_{k+1}
+\shave{\sum^{\left[\!\frac{k}{2}\!\right]}_{i=1}}C^i_{k+1}
\,\bold A_{i+1}\wedge\bold A_{k+1-i}\,\,+\right.\\
&\quad\left.+\,\,\bold A_1\wedge\bold A_{k+1}
+\shave{\sum^{\left[\!\frac{k}{2}\!\right]}_{i=1}}
C^i_{k+1}\,\bold A_i\wedge\bold A_{k+2-i}\right)=0.
\endaligned
\tag6.11
$$
Applying lemma~6.1 to \thetag{6.11}, we derive the following equality for
$d\,\bold A_{k+1}$:
$$
\hskip -2em
\gathered
d\,\bold A_{k+1}=\bold A_0\wedge\bold A_{k+2}+\bold A_1\wedge
\bold A_{k+1}\,+\\
+\sum^{\left[\!\frac{k}{2}\!\right]}_{i=1}
C^i_{k+1}\,\bold A_i\wedge\bold A_{k+2-i}
+\sum^{\left[\!\frac{k+2}{2}\!\right]}_{i=2}
C^{i-1}_{k+1}\,\bold A_i\wedge\bold A_{k+2-i}.
\endgathered
\tag6.12
$$
Comparing \thetag{6.12} and \thetag{6.9}, we can write the following
recurrent formula for $C^i_{k+1}$:
$$
\hskip -2em
C^i_{k+2}=\cases\kern 25pt 1 &\text{for \ }i=0;\\
C^{i-1}_{k+1}+C^i_{k+1} &\text{for \ }0<2i<k+1;\\
\kern 18pt C^{i-1}_{k+1}&\text{for \ }2i=k+1.
\endcases
\tag6.13
$$
Though formula \thetag{6.13} is quite similar to corresponding recurrent
formula for binomial coefficients, it doesn't coincide with that formula.
\par
    Now let's study terms denoted by dots in formula \thetag{6.10}. Total
sum of all these terms is given by the following explicit formula:
$$
\hskip -2em
\gathered
S=\sum^{\left[\!\frac{k}{2}\!\right]}_{i=1}
\sum^{\left[\!\frac{i}{2}\!\right]}_{s=1}C^i_{k+1}
\,C^s_{i+1}\,\bold A_s\wedge\bold A_{i+1-s}\wedge
\bold A_{k+1-i}\,-\\
-\sum^{\left[\!\frac{k}{2}\!\right]}_{r=1}
\sum^{\left[\!\frac{k+1-r}{2}\!\right]}_{e=1}C^r_{k+1}
\,C^e_{k+2-r}\,\bold A_r\wedge\bold A_e\wedge
\bold A_{k+2-r-e}.
\endgathered
\tag6.14
$$
Indices in external product in first sum of formula \thetag{6.14}
satisfy inequalities
$$
\xalignat 2
&\hskip -2em
1\leqslant i<k+1-i,
&&1\leqslant s<i+1-s.
\tag6.15
\endxalignat
$$
Inequalities \thetag{6.15} mean that indices in external product
$\bold A_s\wedge\bold A_{i+1-s}\wedge\bold A_{k+1-i}$ are properly
arranged, i\.\,e\. they are in growing order:
$$
s<i+1-s<k+1-i.
$$
Here are inequalities for indices in external product $\bold A_r
\wedge\bold A_e\wedge\bold A_{k+2-r-e}$:
$$
\xalignat 2
&\hskip -2em
1\leqslant r<k+1-r,
&&1\leqslant e<k+2-r-e.
\tag6.16
\endxalignat
$$
Inequalities \thetag{6.16} cannot provide proper ordering of indices
$r$, $e$, $k+2-r-e$. Therefore we consider three possible subranges
for index $r$:
$$
\align
&\hskip -2em
\text{Subrange 1:}\qquad r<e;
\tag6.17\\
&\hskip -2em
\text{Subrange 2:\qquad} e<r<k+2-r-e;
\tag6.18\\
&\hskip -2em
\text{Subrange 3:\qquad} k+2-r-e<r.
\tag6.19
\endalign
$$
Inequalities \thetag{6.16} define polygon $ABCD$ on $re$-plane
(see Fig\.~6.1 below), sides $AB$ and $AD$ are closed, sides
$BC$ and $CD$ are open. Subranges \thetag{6.17},
\thetag{6.18}, and \thetag{6.19} break this polygon into three
triangular domains $ABE$, $ADE$, and $CDE$.
Segments $AE$ and $DE$ are in open parts of their boundaries.
\par
    {\bf Subrange 1}. In this subrange indices in external product
$\bold A_r\wedge\bold A_e\wedge\bold A_{k+2-r-e}$ are properly ordered.
Therefore we can match them with indices of another external product
$\bold A_s\wedge\bold A_{i+1-s}\wedge\bold A_{k+1-i}$, i\.\,e\. we can
write
$$
\xalignat 3
&\hskip -2em
r=s,
&&e=i+1-s,
&&k+2-r-e=k+1-i.
\quad
\tag6.20
\endxalignat
$$
Third equality in \thetag{6.20} follows from first two ones. Therefore
we can treat first two equalities as a map from $is$-plane to $re$-plane.
This is linear invertible map taking integer points to integer point.
So is inverse map:
$$
\pagebreak 
\xalignat 2
&\hskip -2em
\cases r=s,\\ e=i+1-s,\endcases
&&\cases i=r+e-1,\\s=r.\endcases
\tag6.21
\endxalignat
$$
Due to maps \thetag{6.21} triangle $ABE$ is associated
with triangle $FGH$ (see Fig\.~6.2 below). Indeed, we have the
following correspondence of sides and inequalities:
$$
\xalignat 3
&\hskip -2em
AB\ (r\geqslant 1) &&\longrightarrow &&FG\ (s\geqslant 1);\\
&\hskip -2em
BE\ (e<k+2-r-e) &&\longrightarrow &&GH\ (s>2\,i-k);
\tag6.22\\
&\hskip -2em
EA\ (e>r) &&\longrightarrow &&HF\ (s<i+1-s).\qquad\\
\endxalignat
$$
Due to inequalities in right column of \thetag{6.22} we see that
side $FG$ of subrange 1 mapped to $is$-plane is closed. Other
two sides $GH$ and $HF$ are open.\par
\vskip 0pt\hbox to 0pt{\kern 0pt\hbox{\special{em:graph
pst20a.gif}}\hss}\vskip 210pt
    {\bf Subrange 2}. In this subrange indices in external product
$\bold A_r\wedge\bold A_e\wedge\bold A_{k+2-r-e}$ are not properly
ordered. We need to transpose first two terms in it. As a result we
get external product $\bold A_e\wedge\bold A_r\wedge\bold A_{k+2-r-e}$
that can be matched with external product $\bold A_s\wedge
\bold A_{i+1-s}\wedge\bold A_{k+1-i}$. This yields another pair of
mutually inverse maps linking $re$-plane with $is$-plane. These maps
are given by formulas
$$
\xalignat 2
&\hskip -2em
\cases r=i+1-s,\\ e=s,\endcases
&&\cases i=r+e-1,\\s=e.\endcases
\tag6.23
\endxalignat
$$
Applying \thetag{6.23} to inequalities defining sides of triangle
$AED$, we get
$$
\xalignat 3
&\hskip -2em
AE\ (e<r) &&\longrightarrow &&FH\ (s<i+1-s);\qquad\\
&\hskip -2em
ED\ (r<k+2-r-e) &&\longrightarrow &&HG\ (s>2\,i-k);
\tag6.24\\
&\hskip -2em
DA\ (e\geqslant 1) &&\longrightarrow &&GF\ (s\geqslant 1).\\
\endxalignat
$$
Its important that subrange 2 is mapped onto the same triangle
in $is$-plane as subrange 1, and again side $GF$ is closed,
while other two sides $FH$ and $HG$ of triangle $FHG$ are
open.\par
    {\bf Subrange 3}. In this subrange indices in external product
$\bold A_r\wedge\bold A_e\wedge\bold A_{k+2-r-e}$ also are not properly
ordered. We need to move $\bold A_r$ to third position. Then we get
external product $\bold A_e\wedge\bold A_{k+2-r-e}\wedge\bold A_r$
that can be matched with external product $\bold A_s\wedge\bold A_{i+1-s}
\wedge\bold A_{k+1-i}$. This matching yields two maps inverse to each
other:
$$
\xalignat 2
&\hskip -2em
\cases r=k+1-i,\\ e=s,\endcases
&&\cases i=k+1-r,\\s=e.\endcases
\tag6.25
\endxalignat
$$
Applying \thetag{6.25} to inequalities defining sides of triangle
$DEC$, we get
$$
\xalignat 3
&\hskip -2em
DE\ (k+2-r-e<r) &&\longrightarrow &&GH\ (s>2\,i-k);\\
&\hskip -2em
EC\ (e<k+2-r-e) &&\longrightarrow &&HK\ (s<i+1-s);\qquad
\tag6.26\\
&\hskip -2em
CD\ (r<k+1-r) &&\longrightarrow &&KG\ (2\,i>k+1).\\
\endxalignat
$$
Formulas \thetag{6.26} mean that subrange 3 is mapped onto the
smaller triangle $GHK$ (see Fig\.~6.2). All three sides of this
triangle are open.\par
   Thus, due to \thetag{6.22}, \thetag{6.24}, and \thetag{6.26}
we see that under the action of maps \thetag{6.21}, \thetag{6.23},
and \thetag{6.25} two parts of tetragone $ABCD$ covers triangle
$FGH$ twice, while third part of this tetragone covers smaller
triangle $GHK$. All maps \thetag{6.21}, \thetag{6.23}, and
\thetag{6.25} are given by linear functions with entire coefficients.
Hence they map grid of entire points in $re$-plane onto the grid of
entire points in $is$-plane and vice versa. Note also that inequalities
\thetag{6.15} define triangle $FGK$ complementary to triangle $GHK$
within triangle $FGH$. This means that each entire point of
triangle $FGH$ with closed side $FG$ is associated with three terms
in sum \thetag{6.14}, except for those on segment $GK$. And we have
two terms in sum \thetag{6.14} associated with each inner entire
point of segment $GK$. Therefore in order to prove that $S=0$ in
\thetag{6.14} we should prove series of identities for
coefficients $C^i_k$. First identity
$$
\hskip -2em
C^i_{k+1}\,C^s_{i+1}-C^s_{k+1}\,C^{i+1-s}_{k+2-s}
+C^{i+1-s}_{k+1}\,C^s_{k+1-i+s}=0
\tag6.27
$$
should be fulfilled within open triangle $FGK$. The same identity
\thetag{6.27} should be fulfilled on its side $FG$, except for
ending points $F$ and $G$. Next identity
$$
\hskip -2em
C^{i+1-s}_{k+1}\,C^s_{k+1-i+s}-C^s_{k+1}\,C^{i+1-s}_{k+2-s}
-C^{k+1-i}_{k+1}\,C^s_{i+1}=0
\tag6.28
$$
should be fulfilled within open triangle $GHK$. For exceptional points,
i\.\,e\. for entire points within open segment $GK$, we should prove
the identity
$$
\hskip -2em
C^{i+1-s}_{k+1}\,C^s_{k+1-i+s}-C^s_{k+1}\,C^{i+1-s}_{k+2-s}=0.
\tag6.29
$$
Note that open segment $GK$ has entire points if and only if $k$ is
odd number not less than $7$, i\.\,e\. we should set $k=2\,m+7$, where
$m$ is arbitrary non-negative number. In this case $i=m+4$, while
$s=p+2$, where $p$ is arbitrary non-negative number such that
$2\,p<m+1$. Under these conditions identity \thetag{6.29} reduces to
$$
\hskip -2em
C^{m+3-p}_{2\,m+8}\,C^{p+2}_{m+6+p}
-C^{p+2}_{2\,m+8}\,C^{m+3-p}_{2\,m+7-p}=0.
\tag6.30
$$
In order to prove all these identities we should state formal definition
of coefficients $C^i_k$, other than formula \thetag{6.9}, which is only
a conjecture yet.

\definition{Definition 6.1} Normality coefficients $C^i_k$ are determined
for all integer $k\geqslant 1$ and all integer $i$ such that $0\leqslant
2\,i<k$ by recurrent formula
$$
\hskip -2em
C^i_{k+1}=\cases\kern 25pt 1 &\text{for \ }i=0,\\
C^{i-1}_{k}+C^i_{k} &\text{for \ }0<2i<k,\\
\kern 18pt C^{i-1}_{k}&\text{for \ }2i=k
\endcases
\tag6.31
$$
and by value of initial coefficient $C^0_1=1$ in the series.
\enddefinition
    It is easy to see that definition~6.1 is correct and self-consistent.
Formula \thetag{6.31} is actually the same formula as \thetag{6.13}. Now
let's calculate few initial coefficients in the series and let's arrange
them as a table. Applying formula \thetag{6.31}, we get
$$
\xalignat 6
C^0_1&=1,\\
C^0_2&=1,\\
C^0_3&=1,&C^1_3&=1,\\
C^0_4&=1,&C^1_4&=2,\\
C^0_5&=1,&C^1_5&=3,&C^2_5&=2,\\
C^0_6&=1,&C^1_6&=4,&C^2_6&=5,\\
C^0_7&=1,&C^1_7&=5,&C^2_7&=9,&C^3_7&=5,\\
C^0_8&=1,&C^1_8&=6,&C^2_8&=14,&C^3_8&=14,\\
C^0_9&=1,&C^1_9&=7,&C^2_9&=20,&C^3_9&=28,&C^4_9&=14,\\
C^0_{10}&=1,&C^1_{10}&=8,&C^2_{10}&=27,&C^3_{10}&=48,&C^4_{10}&=42,\\
C^0_{11}&=1,&C^1_{11}&=9,&C^2_{11}&=35,&C^3_{11}&=75,&C^4_{11}&=90,
&C^5_{11}&=42,\\
C^0_{12}&=1,&C^1_{12}&=10,&C^2_{12}&=44,&C^3_{12}&=110,&C^4_{12}&=165,
&C^5_{12}&=132.\\
\endxalignat
$$
One can easily write general formula for elements in first two columns of
this  table:
$$
\xalignat 2
&\hskip -2em
C^0_k=1,
&&C^1_k=k-2.
\tag6.32
\endxalignat
$$
General formula for elements of third column is less obvious:
$$
\hskip -2em
C^2_k=\frac{(k-2)(k-3)}{2}-1.
\tag6.33
$$
However, one can go further and write general formula for all
elements of the table:
$$
\hskip -2em
C^i_k=\prod^i_{s=1}\frac{k-1-s}{s}
-\prod^{k-i}_{s=1}\frac{k-1-s}{s}.
\tag6.34
$$
Formula \thetag{6.34} generalizes \thetag{6.32} and \thetag{6.33}.
In order to prove this general formula it is sufficient to make sure
that it is correct for initial part of the above table and then test
recursion \thetag{6.31} for it. When this is done, proof of the
identities \thetag{6.27}, \thetag{6.28}, and \thetag{6.30} is nothing,
\pagebreak but pure calculations.\par
    Thus, we have proved that $S=0$ in \thetag{6.14}, and hence
we have proved conjecture~6.1. Now we can state it as a theorem.
\proclaim{Theorem 6.1} Differential equations \thetag{5.9} with
$\bold L=\bold A_0$ and $\bold A=\bold A_1$ lead to infinite series
of compatibility conditions \thetag{6.9}, where coefficients $C^i_k$
in \thetag{6.9} are determined by formula \thetag{6.34}.
\endproclaim
\head
7. An example of solution of normality equations.
\endhead
    Theorem~6.1 and formula \thetag{6.9} give a way for constructing
special solutions of normality equations \thetag{5.9}. Let's write
first two equations given by formula \thetag{6.9} and let's loop them
assuming that $\bold A_2=\bold A_0$. Then we have
$$
\xalignat 2
&\hskip -2em
d\,\bold A_0=\bold A_0\wedge\bold A_1,
&&d\,\bold A_1=\bold A_0\wedge\bold A_0=0.
\tag7.1
\endxalignat
$$
Second equation \thetag{7.1} means that $\bold A_1$ is closed $1$-form
within separate fibers of tangent bundle. Locally it is represented
as $\bold A_1=d\varphi$ for some scalar function in $TM$. First equation
\thetag{7.1} for $\bold A_0=\bold L$ then is written as
$$
\hskip -2em
d\,\bold L=\bold L\wedge d\varphi.
\tag7.2
$$
Let's define another $1$-form $\bold M=e^\varphi\,\bold L$. For this
form from \thetag{7.2} we derive:
$$
d\bold M=e^\varphi\,d\varphi\wedge\bold L+e^\varphi\,d\bold L=
e^\varphi\left(d\varphi\wedge\bold L+\bold L\wedge d\varphi\right)=0.
$$
Thus, $\bold M$ appears to be closed form. Like $\bold A_1$ above, it
is determined by some scalar function: $\bold M=dL$. For differential
form $\bold L$ this yields
$$
\hskip -2em
\bold L=e^{-\varphi}\,dL,
\tag7.3
$$
where $\varphi=\varphi(x^1,\ldots,x^n,v^1,\ldots,v^n)$ and
$L=L(x^1,\ldots,x^n,v^1,\ldots,v^n)$. Remember that components
of differential form $\bold L$ determine generalized Legendre
transformation $\lambda$, see \thetag{5.8} and functions
\thetag{1.3}. For these functions from \thetag{7.3} we derive
$$
\hskip -2em
L_i=e^{-\varphi}\,\frac{\partial L}{\partial v^i}.
\tag7.4
$$
{\bf Example}. Let's consider three dimensional case $n=3$ and
let's choose functions
$$
\xalignat 2
&\hskip -2em
\varphi=-v^1,
&&L=v^1+\frac{1}{2}\left((v^2)^2+(v^3)^2\right).
\tag7.5
\endxalignat
$$
Applying formula \thetag{7.4} to functions \thetag{7.5}, in this case
we get
$$
\xalignat 3
&\hskip -2em
L_1=e^{v^1},
&&L_2=v^2\,e^{v^1},
&&L_2=v^3\,e^{v^1}.
\quad
\tag7.6
\endxalignat
$$
These three functions define regular fiber-preserving map \pagebreak
from $TM$ to $T^*\!M$. Its Jacobi matrix can be calculated explicitly.
Indeed, applying \thetag{2.2} to \thetag{7.6}, we get
$$
\hskip -2em
g_{ij}=e^{v^1}\cdot\Vmatrix 1 & 0 & 0\\
v^2 & 1 & 0\\ v^3 & 0 & 1\endVmatrix.
\tag7.7
$$
We also can explicitly calculate inverse matrix for lower triangular
matrix \thetag{7.7}:
$$
\hskip -2em
g^{ij}=e^{-v^1}\cdot\Vmatrix 1 & 0 & 0\\
-v^2 & 1 & 0\\ -v^3 & 0 & 1\endVmatrix.
\tag7.8
$$
Now, using matrix \thetag{7.8}, we apply formula \thetag{2.4} to
components of covector $\bold L$. As a result we get vector $\bold L$
with the following components:
$$
\xalignat 3
&\hskip -2em
L^1=1-(v^2)^2-(v^3)^2,
&&L^2=v^2,
&&L^2=v^3.
\quad
\tag7.9
\endxalignat
$$
Modulus of vector $\bold L$ calculated in non-symmetric metric
\thetag{7.7} is given by formula
$$
\hskip -2em
|\bold L|^2=\sum^n_{i=1}L^i\,L_i=e^{v^1}.
\tag7.10
$$
Now we are able to calculate matrix of projection operator $\bold P$:
$$
P^i_j=\Vmatrix 1-L^1 & -L^1\,v^2 & -L^1\,v^3\\
\vspace{2ex}
-v^2 & 1-(v^2)^3 & -v^2\,v^3\\
\vspace{2ex}
-v^3 & -v^2\,v^3 & 1-(v^3)^2\endVmatrix.
\tag7.11
$$
We used formula \thetag{2.3} for $P^i_j$ and formula \thetag{7.10}
for $|\bold L|^2$. We keep $L^1$ in \thetag{7.11} as notation for
the sake of brevity in order to have formula looking pretty well.
Its value is given by formula \thetag{7.9}.\par
    Next step is to calculate components of tensor field $\bold A$
given by formula \thetag{2.7}. Upon alternating matrix $A^{rs}$
we get the following one:
$$
A^{rs}-A^{sr}=e^{-v^1}\cdot\Vmatrix 0 & v^2 & v^3\\
-v^2 & 0 & 0\\ -v^3 & 0 & 0\endVmatrix.
\tag7.12
$$
Substituting \thetag{7.12} and \thetag{7.11} into \thetag{2.6}, we easily
find that normality equations \thetag{2.6} are fulfilled. Thus, we have
constructed an example of generalized Legendre transformation $\lambda$
satisfying normality equations. It is given by functions \thetag{7.6}.
This is not classical Legendre transformation. However, it differs
from classical one \thetag{1.2} only by scalar factor $e^\varphi$ (see
formula \thetag{7.4}). Therefore we say that \thetag{7.6} is trivial
example of non-classical Legendre transformation satisfying normality
equations.\par
    In order to construct non-trivial solution of normality equations
\thetag{2.6} one should choose another way of looping for the chain of
differential equations \thetag{6.9}. For example we can set $\bold A_3
=\bold A_0$. This leads to more complicated calculations than we carried
out above. Therefore this example will be studied in separate paper.
\head
8. Acknowledgements.
\endhead
     This work is supported by grant from Russian Fund for
Basic Research (project 01-01-00996-a, coordinator of project
Ya\.~T.~Sultanaev), and by grant from Academy of Sciences of
the Republic Bashkortostan (coordinator N.~M.~Asadullin). I am
grateful to these organizations for financial support.

\enddocument
\end